\documentclass[12pt]{amsart}
\usepackage{amssymb}
\usepackage{epsfig}
\usepackage{graphicx}
\numberwithin{equation}{section}

\input xy
\xyoption{all}

\calclayout
\allowdisplaybreaks[3]

\theoremstyle{plain}
\newtheorem{prop}{Proposition}

\newtheorem{theo}[prop]{Theorem}

\newtheorem{lemm}[prop]{Lemma}

\theoremstyle{definition}

\newtheorem{rema}[prop]{Remark}

\def\GL{{\rm GL}}

\def\ra{\rightarrow}

\def\C{{\mathbb C}}

\def\F{{\mathbb F}}

\def\P{{\mathbb P}}
\def\Q{{\mathbb Q}}

\def\Z{{\mathbb Z}}
\def\C{{\mathbb C}}

\def\SL{{\rm SL}}
\def\GL{{\rm GL}}
\def\PGL{{\rm PGL}}

\def\ra{\rightarrow}

\def\PGL{{\rm PGL}}

\def\ra{\rightarrow}

\def\C{{\mathbb C}}
\def\F{{\mathbb F}}
\def\P{{\mathbb P}}
\def\Q{{\mathbb Q}}
\def\R{{\mathbb R}}

\def\Z{{\mathbb Z}}
\def\C{{\mathbb C}}

\makeatother
\makeatletter

\author{Fedor Bogomolov}
\address{Courant Institute\\
                New York University \\
                New York, NY 10012 \\
                USA }
\email{bogomolo@cims.nyu.edu}

\author{Yuri Tschinkel}
\address{Courant Institute\\
                New York University \\
                New York, NY 10012 \\
                USA }
\email{tschinkel@cims.nyu.edu}

\title[Co-fibered products of algebraic curves]{Co-fibered products of algebraic 
curves}

\begin{document}
\date{\today}

\begin{abstract}
We give examples of failure to the existence of co-fibered products in 
the category of algebraic curves.
\end{abstract}

\maketitle
\tableofcontents

\section{Introduction}
\label{sect:introduction}

Let $C_1,C_2$ be smooth complex projective curves. Assume that 
one has a diagram

\centerline{
\xymatrix{
 C \ar[r]^{f_1}\ar[d]_{f_2} & C_1 \\
 C_2   & 
}
} 

\noindent 
with $f_1,f_2$ \'etale surjective morphisms. 
A {\em co-fibered product} is a diagram

\centerline{
\xymatrix{
 C \ar[r]^{f_1}\ar[d]_{f_2} & C_1 \ar[d]^{g_1} \\
 C_2 \ar[r]_{g_2}  & C'
}
} 

\noindent
with $C'$ a curve, $g_1,g_2$ surjective finite morphisms and
$$
g_1f_1=g_2f_2.
$$

The starting point for this note 
was the following question of J. Koll\'ar:
are there obstructions to the existence of 
co-fibered products for unramified covers 
in the category of Riemann surfaces?
In the language of function fields, 
the condition is equivalent to the triviality of the intersection
of the function fields $k(C_1)\cap k(C_2)\subset k(C)$.

More generally, let
$X$ be an algebraic variety over an algebraically closed field $k$
of dimension $n$. Let $K=k(X)$ be its function field.
Consider subfields $k(Y_1),k(Y_2)\subset k(X)$, 
with $\dim(Y_1)= \dim(Y_2)= \dim(X)$. 
We show that under mild conditions on $k$ and the varieties $Y_1,Y_2$, 
one has indeed
$$
k(Y_1)\cap k(Y_2)=k.
$$
This can be achieved as soon as $k^*$ has an element of infinite order.
We will also show that both field extensions 
$k(X)/k(Y_1)$ and  $k(X)/k(Y_2)$
can be unramified, thus satisfying the condition that $f_1$ and $f_2$
above be \'etale.
Using a theorem of Margulis we show 
that co-fibered products for unramified covers
exist unless the curves in question are Shimura curves.
Thus our construction provides {\em all} counterexamples.
In case of curves over the complex numbers, 
we will give examples with small covering degrees, 
e.g., $\deg(g_1)=\deg(g_2)=3$.  

In positive characteristic, related questions on intersections of some specific
function fields have been considered in 
\cite{berkson}, \cite{bremner}, \cite{water}, \cite{beals}, \cite{zieve}. 
A sample result from \cite{berkson} is:
If $k$ is perfect field of characteristic $p$ 
then 
$$
k(x^{pn} + x^{pn-1})\cap k(x^n)\neq k
$$
if and only
if $\gcd(p,n) =1$.
A more precise description of $k(f)\cap k(g)$ is in \cite{bremner}. 
However, the question remains whether co-fibered products exist
for unramified covers of curves over $\bar{\F}_p$.

\section{Elementary examples}

It is easy to construct examples of {\em ramified} covers:
Let $K:=k(\P^1)$ and assume that $k^*$ contains
an element $a$ of infinite order. 
Let us take two involutions
$$
\sigma: x\to 1/x, \quad \text{ and } \quad 
\sigma_a: x \to a/x.
$$ 
They generate a dihedral
group $\mathfrak D_a$ with commutator 
$$
[\sigma, \sigma_2]\, : \,  x\mapsto x/a^2
$$
of {\em infinite} order.
The fields of invariants 
$k(x)^{\sigma}$ and $k(x)^{\sigma_a}$ have index $2$ in $k(x)$,
but the intersection consists of elements which are invariant under 
$\mathfrak D_a$ and hence only of constants.

More generally, let $G\in \PGL_2(k)$ 
be an infinite subgroup generated by two elements of finite
order. The subfields of invariants of $k(\P^1)$ have the 
required property. 
Note that such groups $G$ do not exist for $k$ 
a finite field. Indeed, any finite set of elements in 
$\PGL_2(\bar{\F}_p)$ 
is contained in a subgroup $\PGL_2(\F_q)$, for some $q=p^n$, 
so that this approach fails to produce 
nonintersecting subfields.

\section{Shimura varieties}
\label{sect:general} 

Natural examples of curve covers arise in the theory of arithmetic
groups. Let $G$ be a semi-simple 
algebraic group defined over a number field $F$. 
Fix a model $\mathcal G$ of $G$ over the ring of integers $\mathfrak o_F$.
For every real embedding $\iota\,:\, F\ra \R$ we have a
complex symmetric space 
$$
\mathbb D_{\iota}:=\mathcal G_{\iota}(\R)/{\bf K}_{\iota}, 
\quad \mathbb D:=\prod_{\iota} \mathbb D_{\iota}
$$
where ${\bf K}_{\iota}$ is a maximal compact subgroup.
We have a homomorphism
$$
\phi\,:\, \mathcal G(\mathfrak o_F)\ra \prod_{\iota} \mathcal G_{\iota}(\R). 
$$
Let $\Gamma\subset \phi(\mathcal G(\mathfrak O_K))$ 
be a subgroup of finite index and 
$$
X_{\Gamma}:= \Gamma \backslash \mathbb D.
$$
The quotient $X_{\Gamma}$ is a complex algebraic variety
defined over some finite extension of $\Q$.
For $h\in \phi(\mathcal G(F))$ let 
$\Gamma_h:=h\Gamma h^{-1}$. Then
$\Lambda_h:=\Gamma_h\cap \Gamma$
is a subgroup of finite index in $\Gamma$ and $\Gamma_h$. 
Thus there are surjective maps

\centerline{
\xymatrix{
X:=\mathbb D/\Lambda_h \ar[r]^{f_1}\ar[d]_{f_2} & 
\mathbb D/\Gamma =  X_{\Gamma}\\
\mathbb D/\Gamma_h =X_{\Gamma_h}& 
}
}

\noindent
Both maps are defined over some number field $F'$.
Thus we have two field embeddings
$$
f_1^*(F'(X_{\Gamma}))\subset F'(X), \quad  
f_2^*(F'(X_{\Gamma_h}))\subset F'(X).
$$

\begin{lemm}
\label{lemm:trans}
If $h$ is of infinite order in 
$G(F)$ (modulo the center) then the intersection
$$
f_1^*(F'(X_{\Gamma}))\cap f_2^*(F'(X_{\Gamma_h}))\subset K'(X)
$$ 
is a subfield
of transcendence degree strictly smaller then $\dim(X)$.
If $h$ and $\Gamma$ generate a subgroup of $G(F)$
which acts densely on $\mathbb D/\Gamma$
then the intersection 
$$
f_1^* (F'(X_{\Gamma_h}))\cap f_2^*(F'(X_{\Gamma})) \subset F'(X)
$$
consists only of constants.
\end{lemm}

The same results hold for arbitrary extensions of $F'$, in particular
for complex numbers.

\begin{proof} 
The field $\C(X)$ consists of meromorphic
functions on $\mathbb D:=\prod_\iota D_\iota$ 
which are invariant under the action of $\Gamma_h\cap \Gamma$, and
the subfields 
$$
f_1^*(\C(X_{\Gamma})), f_2^*(\C(X_{\Gamma_h}))
$$ 
of meromorphic functions invariant under
$\Gamma, \Gamma_h$ respectively. 
The intersection 
$f_1^*(\C(X_{\Gamma}))\cap g_2^*(\C(X_{\Gamma_h}))$ consists of functions
invariant under both $\Gamma, \Gamma_h$. If $h$ has infinite
order and if its power is not a central element in $G(F)$ then 
$\Gamma\cap \Gamma_h$ has
infinite index in the group generated by $\Gamma, \Gamma_h$.
This is equivalent to $\C(X_h)$ having an infinite degree over
the intersection $f_1^*(\C(X_{\Gamma}))\cap f_2^*(\C(X_{\Gamma_h}))$. 
Since both fields are algebraic subfields 
the intersection is also algebraic, i.e., 
a finite extension of $\C(y_1,\ldots,y_k)$ for some set $y_1,\ldots,y_k$.
This implies that $k < \dim(X_h)$.

If $\Gamma_h,\Gamma$ 
generate a subgroup which acts on
$\mathbb D$ with a dense orbit
then there are no invariant meromorphic functions on $\mathbb D$.
Thus 
$$
f_1^*(\C(X_{\Gamma}))\cap f_2^*(\C(X_{\Gamma_h}))=\C.
$$
Since the maps $f_i$ are defined over $F'$ the same holds
for arbitrary intermediate subfields $\tilde{F}\subset \C, F'\subset \tilde{F}$.

If the action of $\Gamma$ on $\mathbb D$ is 
cocompact then for a subgroup
of finite index the stabilizers become trivial. 
Then both maps $f_1,f_2$ are finite
unramified covers.
\end{proof}

\section{Curve covers}

Consider $\Gamma:=\SL_2(\Z)$ and $\Gamma_h:= h \,\SL_2(\Z)h^{-1}$,
where $h\in \SL_2(\Q)\setminus \SL_2(\Z)$.
The intersection $\Gamma\cap \Gamma_h$ has finite index in both groups.
However, the action of the group generated by $\Gamma$ and  $\Gamma_h$
on the upper-half plane $\mathfrak H$
is not discrete. In this case we have cusps, i.e., the 
corresponding maps are ramified. A similar argument applies
to any arithmetic group acting on $\mathfrak H$.

Let $D$ be a division algebra of dimension $4$
over $\Q$ which embeds into the $2\times 2$-matrices 
$M_2(\R)$. The spitting field of $D$ is a real-quadratic field
$\Q(\sqrt {d})$.
Let $\Gamma\subset D $ be a subgroup of finite index in the
group of integral quaternions with norm one which does not contain torsion
elements. It acts discretely on 
$\mathfrak H$ with a compact
quotient, the complex points of a projective algebraic curve $C$. 
Let $h\in D$ be an element with a nontrivial denominator and 
$\Gamma_h:=h\Gamma h^{-1}$. 
Write
$$
\Lambda_h:=\Gamma_h\cap \Gamma.
$$
As in Section~\ref{sect:general}, 
we have covers

\centerline{
\xymatrix{
X:=\mathbb D/\Lambda_h \ar[r]^{f_1}\ar[d]_{f_2} & 
\mathbb D/\Gamma =  X_{\Gamma}\\
\mathbb D/\Gamma_h =X_{\Gamma_h}& 
}
}

\noindent
On the other hand,
the group generated by $h$ and $\Gamma$ acts nondiscretely on $\mathfrak H$.
Thus there are no  $h,\Gamma$-invariant elements in the
function field $\C(C')$ and hence no nontrivial common quotient.
The groups $\Gamma$, $\Gamma_h$ and $\Lambda_h$ 
contain no elements of finite order. Hence
they act freely on $\mathbb D$ and the covers $f_1$ and $f_2$ 
are unramified.

\

We now present an series of examples with small covering degrees.
Let $D$ be a quaternion algebra over $\Q$ 
with splitting field $\Q(\sqrt{d})$, for $d> 0$.
Denote by $\Gamma\subset D^1$ the subgroup of 
integer elements in $D$ of norm 1.
Assume that $D$ has the following properties:
\begin{enumerate}
\item $d$ is a square in $\Q_2$ and hence $D\times \Q_2= M_2(\Q_2)$,
\item $\Gamma$ does not contain elements of finite order,
\item $\Gamma$ surjects onto $\SL_2(\Z_2)$.
\end{enumerate}
These conditions are easily satisfied.
Since elements of finite order in $\SL_2(\Z_2)$ have order $2,3,$ and $4$, there are no
elements in $\Gamma$ of order $2$ and $3$, for all  $d > 6$.
Note that $D$ is dense in $\GL_2(\Q_2)$. Let $h\in D$ be an element
which modulo 4 is equal to
$$
\left(\begin{array}{cc} 1 &  -1/2 \\ 
                        0 &  1\end{array}\right),
$$
and put $\Gamma_h:=h\Gamma h^{-1}$. 
The intersection 
$$
\Lambda_h:= \Gamma_h\cap \Gamma
$$ 
contains a subgroup $x= 1 \mod 2$
and a subgroup modulo 4 generated by 
$$
\left(
\begin{array}{cc}
1  & 1 \\
0 & 1 \end{array}\right).
$$
This is a subgroup of index $3$ in $\Gamma$ and also in $\Gamma_h$.
The group $\Lambda_h$ is the preimage of a congruence subgroup in $\SL_2(\Z_2)$.
Since $\SL_2(\Z/2) = \mathfrak S_3$ 
and the unipotent subgroup is $\Z_2$ we obtain
that $\Lambda_h$ has index $3$.

The construction shows that 
it suffices to assume that $\Gamma $ is 
a subgroup of finite index in the group of
integral elements in $D$ of norm 1, 
which surjects onto $\SL_2(\Z_2)$ and has no
torsion. 
For example, we can take any congruence subgroup, i.e., 
insist that $g\in \Gamma$ satisfies $g = 1\mod p$, 
for some prime $p\neq 2$.

\begin{rema}
We cannot achieve that both $g_1,g_2$ are of degree $2$
and unramified.
Indeed, in this case the corresponding extensions 
would be Galois, and the actions of both 
$\Z/2$ could be realized inside an action of a finite group $H$ on $C$.
Thus $k(C')\cap k(C'')$ contains
$k(C)^H$, a nontrivial field.
\end{rema}

Let  $G$ be a semi-simple algebraic group over $\Q$ and
$$
{\rm Comm}_G(\Gamma):=\left\{\, g\in G(\R)\,|\, [\Gamma: \left (g\Gamma g^{-1}\cap \Gamma\right)] <\infty\, \right\}.  
$$
This is a well-defined subgroup of $G(\R)$ containing $\Gamma$. 

Assume that $\Gamma\subset \SL_2(\R)$ is a discrete cocompact subgroup without
torsion elements. Then $X$ admits maps $f_1,f_2$
as above, if and only if $\Gamma$ is an
arithmetic subgroup of $\SL_2(\R)$. 
Note that
$$
[{\rm Comm}_{\SL_2(\R)}(\Gamma) :\Gamma]<\infty
$$ 
{\em unless} $\Gamma$ is arithmetic. Indeed, we have the following 

\begin{theo}\cite[Theorem (B), p. 298]{Margulis-book}
Let $G$ be a semi-simple group over $\mathbb Q$ and 
$\Gamma\subset G(\R)$ an irreducible lattice.
Assume that $\Gamma$
\begin{enumerate}
\item[(i)] is of infinite index in ${\rm Comm}_G(\Gamma)$;
\item[(ii)] is finitely generated;
\item[(iii)] satisfies property $QD$.
\end{enumerate}
Then $\Gamma$ is arithmetic.
\end{theo}

In our applications, $\Gamma$ automatically satisfies properties (ii) and (iii)
(see \cite[Chapter IX]{Margulis-book} for definitions and results).

\bibliographystyle{smfalpha}
\bibliography{shimura}

\end{document}